\documentstyle{article}

\title{Mathematics Journals Should Be Electronic and Free
\thanks{Redistributed from
\emph{Notices Amer. Math. Soc.}, 44 (1997), no. 8, 2}}
\author{Steve Krantz}
\date{}
\begin{document}
\maketitle

\thispagestyle{empty}

Of course \emph{I} don't believe what the title says. But I got your
attention.

Journals record what we know and preserve our collected knowledge.  They
establish who proved what. A dean checks journals to verify that a tenure
candidate publishes in established, archived forums of good repute.
\emph{All of these journal functions would be lost if we were
immediately to replace all paper journals with electronic media}.
Are you shocked? Then take heed.

\begin{quote}
Nobody knows how to archive electronic media. The expected time value
of data stored on a CD-ROM is ten years. Tapes oxidize and decompose.
Hard discs crash. Zip cartridges and Jaz cartridges and Bernoulli
cartridges are only warranted for one to five years. An individual
easily can keep his few megs of data alive and on current media; a
library would be overwhelmed by the task.

The hardware protocol changes every few years. (If you found an Edison
cylinder labeled ``Proof of the Riemann Hypothesis'', you could read
the label, but how would you access the proof?)

Imagine presenting a tenure case to your dean (a philologist of
Tlingit), and that the candidate's work is all published in electronic
math journals. You assert that (i) these are scholarly journals with
distinguished editors, (ii) they are refereed and archived. Good luck.
Over many years you might educate the dean. But you might not. (You
know how difficult it is to get a tenure candidate past even the
department if most of his publications are in ``unrefereed''
conference proceedings; if the publications are in ``free'' electronic
journals, then exponentiate that difficulty.)

In one hundred years there will be no CD-ROMs, no \TeX, no DOS, no
Windows95, and no PCs. Anything stored and archived today will be
inaccessible then. Some say, ``We'll plan to change media and software
regularly.'' Nice thought, but there are few such mechanisms in place;
when money gets tight, corners will be cut, and vast amounts of the
literature will be lost.
\end{quote}

You can always pull a book or bound journal off the shelf and read it. You need
only a pair of eyes and an education (both time-tested devices). The archival
value of printed books, especially those printed on acid-free paper, is well
established.

People hear ``electronic journal'' and think ``free journal''. Nonsense. Most
of the cost of producing a journal comes from typesetting, formatting,
editorial work, clerical work, accessibility, and archiving. Electronic media
may save \emph{some} costs, but they will not
change the landscape.

Some argue that mathematicians would be willing to give up a great deal in the
\emph{quality} of journals if such a sacrifice
would substantially reduce the price.  I find this curious. I do not see many
mathematicians driving Yugos, or wearing hopsacks, or doing their e-mail and
compiling \TeX\ on a Radio Shack TRS-80 just because it is cheaper. I do not
think we are accurately assessing our own values. If you proved Fermat's Last
Theorem, would you want it typed up on a vintage 1915 Underwood, with the
mathematical symbols written in by hand, duplicated with a ditto master, spiral
bound, and dropped in bundles from low-flying planes? I think it is time that
we think carefully about what we expect from a journal.

We are slowly being co-opted by electronic media: (i) our papers used to be
typed \emph{for} us, but now we do it ourselves;
(ii) many of us volunteer time to help maintain the departmental computer
system; (iii) mathematicians now give (considerable) time to running so-called
``free'' electronic journals. Perhaps the benefits of electronic media are
offset by their negative impact on the infrastructure of the profession.

Electronic media disseminate information quickly, and that is
important. As for electronic media being the basis for a publishing
revolution? Not yet. The emperor is still dressing.

\end{document}